\documentclass[12pt]{amsart}
\usepackage[latin1]{inputenc}

\usepackage{pstricks}

\usepackage{ifthen}
\newcommand{\mymarginpar}[1]{%
   \marginpar{\ifthenelse{\isodd{\arabic{page}}}{\flushleft #1}{\flushright #1}}}
\newcommand{\QQQ}[1]{ \mymarginpar{\textbf{QQQ}} #1 }
\renewcommand{\QQQ}[1]{}
\usepackage{longtable}

\usepackage{enumerate}

\theoremstyle{plain} 
\newtheorem{Theorem}{Theorem}[section]

\newtheorem{Proposition}[Theorem]{Proposition}

\theoremstyle{definition} 
\newtheorem{Definition}[Theorem]{Definition}
\newtheorem{Remark}[Theorem]{Remark}

\numberwithin{equation}{section}

\newcommand{\Hankel}{\mathcal{H}}
\newcommand{\Lukas}{\mathcal{L}}
\newcommand{\Lukasirr}{\Lukas^{\mathrm{irr}}}
\newcommand{\Motzkin}{\mathcal{M}}
\newcommand{\Motzkinirr}{\Motzkin^{\mathrm{irr}}}
\newcommand{\abs}[1]{\left\lvert #1 \right\rvert}
\DeclareMathOperator{\sign}{{\mathrm sign}}                     
\newcommand{\SG}{\mathfrak{S}}                   
\newcommand{\IN}{\mathbf{N}}                     
\newcommand{\IZ}{\mathbf{Z}}

\newcommand{\alg}[1]{\mathcal{#1}}                
\renewcommand{\phi}{\varphi}

\begin{document}

\title{Cumulants, Lattice Paths, and Orthogonal Polynomials}
\author{Franz Lehner}


\address{
Franz Lehner\\
In\-sti\-tut f\"ur Mathe\-ma\-tik C\\
Tech\-ni\-sche Uni\-ver\-si\-t\"at Graz\\
Stey\-rer\-gas\-se 30, A-8010 Graz\\
Austria}
\email{lehner@finanz.math.tu-graz.ac.at}

\date{\today}

\begin{abstract}
A formula expressing free cumulants in terms of Jacobi parameters of
the corresponding orthogonal polynomials is derived.
It combines Flajolet's theory of continued fractions and
the Lagrange inversion formula.
For the converse we discuss Gessel-Viennot theory to express Hankel determinants in
terms of various cumulants.
\end{abstract}

\maketitle
\section{Introduction}

In his seminal paper \cite{Flajolet:1980:continuedfractions} Flajolet
gave a combinatorial interpretation of continued fractions and derived
generating function expressions for enumeration problems of various
lattice paths.
This is connected to work of Karlin and McGregor
and others on birth and death processes, see 
\cite{FlajoletGuillemin:2000:processes} or
\cite{IsmailLetessierMassonValent:1990:processes} for a survey.
The basic principle has been rediscovered many times, see
e.g.\ \cite{Gerl:1984:continuedfraction} for an application to random walks
and \cite{AccardiBozajko:1998:Gaussianization} for an interpretation
in noncommutative probability.
Since Stieltjes' times continued fractions have also been a basic ingredient
of the theory of orthogonal polynomials.
A synthesis of both aspects can be found in Viennot's memoir 
\cite{Viennot:1983:polynomes} or in the survey \cite{Viennot:1985:polynomials}.

Orthogonal polynomials and their Jacobi operators have been extensively
studied in connection with moment problems and spectral theory
\cite{Akhiezer:1965:momentproblem}. 
Spectral theory of convolution operators was also one of Voiculescu's
movitations for the development of free probability theory
\cite{VoiculescuDykemaNica:1992:free} 
and it would be interesting to understand the behaviour of orthogonal
polynomials under free convolution.
The special case of free projections, which is up to a translation
equivalent to the case of free generators,
have been considered in the literature, see e.g.\ 
\cite{Promislow:1978:dimension,
CohenTrenholme:1984:orthogonal,
Trenholme:1988:Green,
AkiyamaYoshida:1999:orthogonal}.

In this note we derive a formula for free cumulants in the spirit
Flajolet's and Viennot's theory of lattice paths.
After discussions with A.~Lascoux, he gave a proof in terms of symmetric
functions of a more general formula for
the coefficients of powers of continued fractions, cf.\ 
\cite{Lascoux:2000:motzkin,Zeng:2000:powers}. 

The paper is organized as follows.

In section~\ref{sec:LatticePaths} we review Flajolet's formula for
the generating function of Motzkin paths.

In section~\ref{sec:orthogonalpolynomials} we briefly discuss the relevant
facts about orthogonal polynomials.

In section~\ref{sec:cumulants} we review the definitions of
various cumulants.

In section \ref{sec:freecumulantsformula} we prove the main result,
a formula expressing free cumulants in terms of Jacobi parameters.

Finally in section~\ref{sec:Hankeldeterminants} we indicate how Gessel-Viennot
theory of Hankel determinants can be used to express the Jacobi parameters in
terms of cumulants.

\section{Enumeration of Lattice Paths}
\label{sec:LatticePaths}

\begin{Definition}
  A \emph{lattice path} is a sequence of points in the integer lattice $\IZ^2$.
  A pair of consecutive points is called a \emph{step} of the path.
  A \emph{valuation} is a function on the set of possible steps
  $\IZ^2\times\IZ^2$.
  A valuation of a path is the product of the valuations of its steps.
  In the rest of this paper all lattice paths will have the following properties.
  \begin{enumerate}
   \item Starting point and end point lie on the $x$-axis.
   \item The $y$-coordinates of all points are nonnegative.
   \item In each step, the $x$-coordinate is incremented by one.
  \end{enumerate}
  Thus a path of length $n$ will start at some point $(x_0,0)$
  and end at $(x_0+n,0)$.
  The valuations will be independent of the $x$-coordinates of the points.
  Therefore the $x$-coordinates are redundant and we can
  represent a path $\pi$ by the sequence of its $y$-coordinates
  $(\pi(0),\pi(1),\dots,\pi(n))$.
  We call a path \emph{irreducible} if it does not touch the
  $x$-axis except at the start and at the end.
  Every path has a unique factorization into irreducible ones.
\end{Definition}
We will be concerned with two types of paths.
\begin{Definition}
  \begin{enumerate}[(i)]
   \item 
    A \emph{Motzkin path} of length $n$ is a lattice path starting at $P_0 = (0,0)$ and
    ending at $(n,0)$, all of whose $y$ coordinates are nonnegative and
    whose steps are of the three following types.
    \begin{itemize}
     \item [] \emph{rising step}: $(1,1)$ 
     \item [] \emph{horizontal step}: $(1,0)$ 
     \item [] \emph{falling step}: $(1,-1)$ 
    \end{itemize}
    Motzkin paths without horizontal steps are called \emph{Dyck paths}.
    We will denote the set of Motzkin paths by $\Motzkin_n$
    and the subset of irreducible Motzkin paths by $\Motzkinirr_n$.
    Motzkin paths are counted by the well known \emph{Motzkin numbers}.
   \item A \emph{\L{}ukasiewicz path} of length $n$ is a path starting at $(0,0)$
    and ending at $(n,0)$ whose steps are of the following types.
    \begin{itemize}
     \item [] \emph{rising step}: $(1,1)$ 
     \item [] \emph{horizontal step}: $(1,0)$ 
     \item [] \emph{falling steps}: $(1,-k)$, $k>0$ 
    \end{itemize}
    We denote the set of \L{}ukasiewicz paths of length $n$ by $\Lukas_n$
    and the subset of irreducible \L{}ukasiewicz paths by $\Lukasirr_n$.
    \L{}ukasiewicz paths form a \emph{Catalan family}, i.e.,
    they are counted by the \emph{Catalan numbers}.
  \end{enumerate}
\end{Definition}
The name \L{}ukasiewicz path is motivated by the natural bijection to the 
\emph{\L{}ukasiewicz language}, the language understood by
calculators using reverse polish notation,
which was used by Raney in \cite{Raney:1960:composition}
to give a combinatorial proof of the Lagrange inversion formula.

Flajolet's formula expresses the generating function of weighted Motzkin paths
as a continued fraction.
\begin{Theorem}[{\cite{Flajolet:1980:continuedfractions}}]
  \label{thm:Flajolet}
  Let
  \begin{equation}
    \label{eq:Flajoletformula}
    \mu_n = \sum_{\pi\in \Motzkin_n}
             v(\pi)
  \end{equation}
where the sum is over the set of Motzkin paths
$\pi=(\pi(0)\dots\pi(n))$ of length $n$.
Here $\pi(j)$ is the level after the $j$-th step,
and the \emph{valuation} of a path is the product
of the valuations of its steps $v(\pi)=\prod_1^n v_i$,
the latter being
\begin{equation}
  \label{eq:Flajoletvaluation}
  v_i = v(\pi(i-1),\pi(i)) =
  \begin{cases}
    1                  &\text{if the $i$-th step rises}\\
    a_{\pi(i-1)}       &\text{if the $i$-th step is horizontal}\\
    \lambda_{\pi(i-1)} &\text{if the $i$-th step falls}
  \end{cases}
\end{equation}
Then the generating function
$$
M(z) = \sum_{n=0}^\infty \mu_n z^n
$$
has the continued fraction expansion
\begin{equation}
  \label{eq:Flajolet:contfrac}
  M(z) = \cfrac{1}%
           {1 - \alpha_0 z - \cfrac{\lambda_1 z^2}%
                               {1 - \alpha_1 z-\cfrac{\lambda_2z^2}%
                                                 {\ddots}}}
\end{equation}
\end{Theorem}

\section{Orthogonal Polynomials}
\label{sec:orthogonalpolynomials}

A sequence of (formal) orthogonal polynomials
is a sequence of monic polynomials $P_n(x)$ of degree $\deg P_n = n$
together with a linear functional $\mu$ on the space of polynomials
with moments $\mu(x^n)=\mu_n$ such that
$\mu(P_m(x)\,P_n(x)) = \delta_{mn} s_n$ for some coefficients $s_n$.
We will always assume that $s_0=1$, that is $\mu(x^0)=1$.
It is then easy to see that such a sequence satisfies a three-term 
recurrence relation
$$
xP_n(x) = P_{n+1}(x) + a_{n}P_n(x) + \lambda_n P_{n-1}(x)
$$
with $\lambda_0=0$. The numbers $a_n$, $\lambda_n$ are the so-called
\emph{Jacobi parameters}.
One verifies easily that the polynomials
$P_n(x)=\frac{D_n(x)}{\Delta_{n-1}}$ satisfy the orthogonality condition,
where
\begin{equation}
  \label{eq:Dn(x)}
  D_n(x)
  = \begin{vmatrix}
      \mu_0    & \mu_1 & \dots & \mu_n \\
      \mu_1    & \mu_2 &       & \mu_{n+1} \\
      \vdots   &       &       & \vdots \\
      \mu_{n-1}&       & \dots & \mu_{2n-1} \\
      1        & x     & \dots & x^n
    \end{vmatrix}
\end{equation}
and
\begin{equation}
  \label{eq:Deltan(x)}
  \Delta_n
  = \begin{vmatrix}
      \mu_0    & \mu_1    & \dots & \mu_n \\
      \mu_1    & \mu_2    &       & \mu_{n+1} \\
      \vdots   &          &       & \vdots \\
      \mu_n    & \mu_{n+1}& \dots & \mu_{2n}
    \end{vmatrix}
\end{equation}
is the $n$-th \emph{Hankel determinant}.

Coming from an operator theoretic background, the most natural way to
make the connection to Flajolet's formula is perhaps via matrices.
The Jacobi matrix model for the moment functional is
$$
J = 
\begin{bmatrix}
  \alpha_0 & 1 \\
  \lambda_1& \alpha_1 & 1\\
           & \lambda_2& \alpha_2 & 1 \\
           &          & \ddots   & \ddots    &\ddots
\end{bmatrix}
,
$$
that is, if we denote the basis of the vector space by $e_0, e_1,\dots$,
with inner product $\langle e_m, e_n \rangle = \delta_{mn} s_n$,
then $Je_n = e_{n+1} + \alpha_n e_n + \lambda_n e_{n-1}$.
Then it is easy to see that $P_n(J)\,e_0 = e_n$ and therefore $J$
satisfies $\langle J^n e_0,e_0 \rangle = \mu_n$.
Expanding the matrix power yields
$$
\mu_n = \sum_{i_1,\dots,i_{n-1}\ge 0}
         J_{0,i_1} J_{i_1,i_2} \dots J_{i_{n-1},0}
$$
and because of tridiagonality the sum is restricted to indices
$|i_{j+1}-i_j|\le 1$. The summands
$J_{0,i_1} J_{i_1,i_2} \dots J_{i_{n-1},0}$
can be interpreted as the
valuations of the Motzkin paths $(0,i_1,i_2,\dots,i_{n-1},0)$
with weights $v(y,y')=J_{y+1,y'+1}$, where
 $J_{i,i-1}=\lambda_i$, $J_{i,i}=\alpha_i$, $J_{i,i+1}=1$
and this is exactly Flajolet's formula \eqref{eq:Flajoletformula}.

\section{Cumulants}
\label{sec:cumulants}

Cumulants linearize convolution of probability measures coming from
various notions of independence.
\begin{Definition}[{\cite{VoiculescuDykemaNica:1992:free}}]
  A \emph{non-commutative probability space} is pair $(\alg{A},\phi)$
  of a (complex) unital algebra $\alg{A}$ and a unital linear functional $\phi$.
  The elements of $\alg{A}$ are called \emph{(non-commutative) random variables}.
  The collection of moments $\mu_n(a) = \phi(a^n)$
  of such a random variable $a\in\alg{A}$ will be called
  its \emph{distribution} and denoted $\mu_a = (\mu_n(a))_n$.
\end{Definition}
Thus noncommutative probability follows the general ``quantum'' philosophy
of replacing function algebras by noncommutative algebras.
We will review several notions of independence below.
Convolution is defined as follows. Let $a$ and $b$ be ``independent'' random variables.
Then the convolution of the distributions of $a$ and $b$ is defined to be
the distribution of the sum $a+b$. In all the examples below, the distribution
of the sum of  ``independent'' random variables only depends on the individual
distributions of the summands and therefore convolution is well defined
and the $n$-th moment $\mu_n(a+b)$ is a polynomial function of the moments
of $a$ and $b$ of order less or equal to $n$.
For our purposes it is sufficient to axiomatize cumulants as follows.
\begin{Definition}
  \label{def:cumulants}
  Given a notion of independence on a noncommutative probability space $(\alg{A},\phi)$,
  a sequence of maps $a\mapsto k_n(a)$, $n=1,2,\dots$ is called a \emph{cumulant sequence} if 
  it satisfies
  \begin{enumerate}
   \item additivity: if $a$ and $b$ are ``independent'' random variables,
    then $k_n(a+b)=k_n(a)+k_n(b)$.
   \item homogenity: $k_n(\lambda a) = \lambda^n k_n(a)$.
   \item $k_n(a)$ is a polynomial in the first $n$ moments of $a$ with leading term $\mu_n(a)$.
    This ensures that conversely the moments can be recovered from the cumulants.
  \end{enumerate}
\end{Definition}

We will review here \emph{free}, \emph{classical} and \emph{boolean}
cumulants via their matrix models.
The reader interested in ``$q$-analogues'' is referred to
the $q$-Toeplitz matrix models of \cite{Nica:1995:oneparameter} and 
\cite{Nica:1996:crossingsembracings} which yield similar
valuations on \L{}ukasiewicz paths.

\subsection{Free Cumulants}

Free probability was introduced by Voiculescu in~\cite{Voiculescu:1985:operator}
and has seen rapid development since, see~\cite{VoiculescuDykemaNica:1992:free} and
 the more recent survey~\cite{Voiculescu:2000:lectures}.
\begin{Definition}
  Given a noncommutative probability space $(\alg{A},\phi)$, the
  subalgebras $\alg{A}_i\subseteq \alg{A}$ are called
  \emph{free independent} (or \emph{free} for short) if
  \begin{equation}
    \label{eq:FreenessDefinition}
    \phi(a_1a_2\cdots a_n)=0
  \end{equation}
  whenever $a_j\in \alg{A}_{i_j}$ with
  $\phi(a_j)=0$ and $i_j\ne i_{j+1}$ for $j=1,\dots,n-1$.
  Elements $a_i\in\alg{A}$ are called free if the unital subalgebras generated 
  by $a_i$ are free.
\end{Definition}

Existence of free cumulants with the properties of definition~\ref{def:cumulants}
was proved already in~\cite{Voiculescu:1985:operator}.
A beautiful systematic theory of free cumulants was found by Speicher in his
combinatorial approach to free probability via non-crossing partitions
\cite{Speicher:1994:multiplicative}; see \cite{Lehner:2001:cumulants}
for an explanation why noncrossing partitions appear.
The explicit computation involves generating functions as follows.
\begin{Theorem}[{\cite{Voiculescu:1986:addition}}]
  Let $M(z)=1+\sum_{n=1}^\infty \mu_n z^n$ be the ordinary moment generating function
  and $C(z)=1+\sum_{n=1}^\infty c_n z^n$ be the function 
  implicitly defined by the relation $C(zM(z))=M(z)$.
  Then the coefficients $c_n$ satisfy the requirements of definition~\ref{def:cumulants}
  and are called the \emph{free} or \emph{non-crossing cumulants}
  and $C(z)$ is the \emph{cumulant generating function}.
\end{Theorem}

The moments and the free cumulants are related explicitly by the following combinatorial
formula 
(see \cite{Raney:1960:composition,Cori:1983:wordstrees} for
the formulation in terms of \L{}ukasiewicz language, and
\cite{Speicher:1994:multiplicative} for non-crossing partitions).
Define a valuation on \L{}ukasiewicz paths by putting the following
weights on the steps
\begin{equation}
  \label{eq:Lukasvaluation:NC}
  \begin{aligned}
    v(y,y-k) &= c_{k+1} \qquad k\geq 0 \\
    v(y,y+1) &= 1,
  \end{aligned}
\end{equation}
then
\begin{equation}
  \label{eq:mun=lukasiewicz}
  \mu_n = \sum_{\pi\in\Lukas_n} v(\pi)
  .
\end{equation}

Another interpretation of this relation is the \emph{Fock space model}
of Voiculescu \cite{Voiculescu:1986:addition}
(see \cite{Haagerup:1997:RStransforms} for a simpler proof).
It involves Toeplitz operators as follows.
Let $S$ be the forward shift on $\ell_2(\IN_0)$, i.e., $Se_n=e_{n+1}$.
A (formal) \emph{Toeplitz operator} is a linear combination of powers
of $S$ and of its adjoint, the backward shift $S^* e_n = e_{n-1}$.
The linear functional $\omega(X) = \langle X e_0, e_0 \rangle$ is
called the \emph{vacuum expectation}.
Let $c_n$ be the free cumulants of the moment sequence $(\mu_n)$
as defined above and set
\begin{equation}
  \label{eq:ToeplitzModel}
  T=S^*+\sum_{n=0}^\infty c_{n+1} S^{n}
  .
\end{equation}
Then $\omega(T^n) = \mu_n$ 
and writing $T$ in matrix form
$$
T = 
\begin{bmatrix}
  c_1   & 1  \\ 
  c_2   & c_1 & 1 \\
  c_3   & c_2 & c_1    & 1 \\
  \vdots&     &\ddots&   &\ddots
\end{bmatrix}
$$
we can expand the matrix product and obtain
$$
\omega(T^n)
     = \langle T^n e_0,e_0 \rangle
     = \sum_{i_1,\dots,i_{n-1}\ge 0}
        T_{0,i_1} T_{i_1,i_2} \dots T_{i_{n-1},0}
$$
where $T_{ij}\ne 0$ only for $j\ge i-1$.
Again this can be interpreted as a sum over lattice paths
which this time turn out to be \L{}ukasiewicz paths
with weights \eqref{eq:Lukasvaluation:NC} coming from the matrix entries
$$
T_{i,i+1}=1 \qquad T_{i,i-k}=c_{k+1},\ k\geq 0
$$
and we obtain the sum \eqref{eq:mun=lukasiewicz}.
There are multivariate generalizations of this Toeplitz model,
see~\cite{Nica:1996:Rtransforms}.

\subsection{Classical Cumulants}

Classical cumulants linearize convolution of measures and can be
defined via the Fourier transform.
\begin{Definition}
  Let $F(z) = \mu(e^{xz}) = \sum_{n=0}^\infty \frac{\mu_n}{n!} z^n$ be
  the exponential moment generating function.
  Let $K(z) = \log F(z) = \sum_{n=1}^\infty \frac{\kappa_n}{n!} z^n$
  be its formal logarithm.
  The coefficients $\kappa_n$ are the \emph{classical cumulants} of
  the functional $\mu$.
\end{Definition}

There is a model on bosonic Fock space for classical cumulants
which is analogous to the Toeplitz model \eqref{eq:ToeplitzModel}.

Let $D$, $x$ be annihilation and creation operators which satisfy
the \emph{canonical commutation relations} ($CCR$):
$[D,x]=1$, i.e., 
on the Hilbert space with basis $e_n$ and
inner product $\langle e_m, e_n \rangle = n! \delta_{m,n}$,
set
$D e_n = n e_{n-1}$, $x e_n = e_{n+1}$.
Then $D$ and $x$ are adjoints of each other.

Denote again by $\omega$ the vacuum expectation 
$\omega(T)=\langle T e_0, e_0 \rangle$, 
i.e.\ $\omega(x^k D^n)=0$ $\forall n\ne0$
and $\omega(f(x))=f(0)$.
Then the Fourier-Laplace transform of 
$$
T=D+\sum_{n=0}^\infty \frac{\kappa_{n+1}}{n!} x^n
$$
is 
$$
\omega(e^{zT})=e^{\sum_{n=1}^\infty \frac{\kappa_n}{n!} z^n}.
$$
In the basis $\{e_n\}$, $T$ has matrix representation
$$
\hat{T} = 
\begin{bmatrix}
  \kappa_1           & 1    \\
  \kappa_2           & \kappa_1          & 2  \\
  \frac{\kappa_3}{2!}& \kappa_2          & \kappa_1 & 3\\
  \frac{\kappa_4}{3!}&\frac{\kappa_3}{2!}& \kappa_2 & \kappa_1 & 4  \\
  \vdots             &                   &          & \ddots   &   &\ddots
\end{bmatrix}
$$
e.g.\ for the stnadard Gaussian distribution we have $\kappa_2=1$
and higher cumulants vanish, so the model is tridiagonal
$$
\begin{bmatrix}
  0 & 1 \\
  1 & 0 & 2 \\
    & 1 & 0    & 3 \\
    &   &\ddots&\ddots&\ddots
\end{bmatrix}
$$
and this is the well known Jacobi operator matrix for the Hermite polynomials.
Similar to the formula above we get the following combinatorial sum
for the moments of $T$:
$$
\mu_n = \sum_{\pi\in\Lukas_n} v(\pi)
$$
with valuation
\begin{equation}
  \label{eq:Lukasvaluation:classical}
  \begin{aligned}
    v(y, y-k) &= \frac{\kappa_{k+1}}{k!} \qquad k\geq 0 \\
    v(y, y+1) &= y+1
  \end{aligned}
\end{equation}

\QQQ{Compare this to Flajolet's and Biane's bijections between Motzkin paths and partitions}

\subsection{Boolean Cumulants}

Boolean cumulants linearize \emph{Boolean convolution} which was introduced in
\cite{SpeicherWoroudi:1997:boolean} 
(compare also \cite{vonWaldenfels:1973:approach}, \cite{vonWaldenfels:1975:intervalpartitions})
and is a special case of 
\emph{conditional independence}
(\cite{BozejkoSpeicher:1991:psiindependent},
\cite{BozejkoLeinertSpeicher:1996:convolution}).
In the theory of random walks boolean cumulants arise as
\emph{first return probabilities}.
\begin{Definition}
  Given a noncommutative probability space $(\alg{A},\phi)$, the
  subalgebras $\alg{A}_i\subseteq \alg{A}$ are called
  \emph{boolean independent} if
  \begin{equation}
    \label{eq:BooleanDefinition}
    \phi(a_1a_2\cdots a_n) = \phi(a_1)\, \phi(a_2) \cdots \phi(a_n)
  \end{equation}
  whenever $a_j\in \alg{A}_{i_j}$ with
  $i_j\ne i_{j+1}$ for $j=1,\dots,n-1$.
\end{Definition}
With this notion of independence, convolution of measures is well defined and
the appropriate cumulants can be calculated as follows
\cite{SpeicherWoroudi:1997:boolean}.
Let $M(z) = \sum_{n=0}^\infty \mu_n z^n$ be the ordinary moment generating
function. Then
\begin{equation}
  \label{eq:booleancumulantsdefinition}
  M(z) = \frac{1}{1-H(z)}  
\end{equation}
where 
\begin{equation}
  \label{eq:booleancumulants}
  H(z) = \sum_{n=1}^\infty h_n z^n
  ,
\end{equation}
where $h_n$ are the
\emph{boolean cumulants}
of the distribution with moment sequence $(\mu_n)$.
Expanding \eqref{eq:booleancumulantsdefinition} in a geometric series
we find that the moments can be expressed as Cauchy convolution
$$
\mu_n = \sum_{r=1}^n
        \sum_{\substack{i_1+\dots+i_r=n\\ i_j\geq 1} }
         h_{i_1} h_{i_2} \cdots h_{i_r}
$$
and from Theorem~\ref{thm:Flajolet} and the fact that every Motzkin path has
a unique factorization into irreducible paths, it follows that
\begin{equation}
  \label{eq:boolean=summotzkinirr}
  h_n = \sum_{\pi\in \Motzkinirr_n} v(\pi)  
\end{equation}
where $v$ is the valuation \eqref{eq:Flajoletvaluation}.
This can also be seen by comparing the continued fraction expansion
\eqref{eq:Flajolet:contfrac} with
\eqref{eq:booleancumulantsdefinition}.
In terms of the Jacobi parameters the boolean cumulant generating
function has the continued fraction expansion
$$
H(z) = \alpha_0 z + \cfrac{\lambda_1 z^2}%
                      {1-\alpha_1 z - \cfrac{\lambda_2 z^2}%
                                        {1-\alpha_2 z
                                          - \cfrac{\lambda_3 z^2}{\ddots}}}
.
$$

Similarly, we can express the boolean cumulants in terms of free cumulants
(resp.\ classical cumulants)
\begin{equation}
  \label{eq:boolean=sumLukasirr}
  h_n = \sum_{\pi\in \Lukasirr_n} v(\pi)  
\end{equation}
using the valuations
\eqref{eq:Lukasvaluation:NC}
(resp.\ \eqref{eq:Lukasvaluation:classical}).

\section{A Formula for Free Cumulants}
\label{sec:freecumulantsformula}

\begin{Theorem}
  For $n\geq2$, we have the following formula for the free cumulant $c_n$
  in terms of the Jacobi parameters:
  \begin{equation}
    \label{eq:freecumulant=motzkin}
    c_n = \sum_{\pi\in\Motzkin_n}
          \frac{(-1)^{\abs{\pi}_0-1}}{n-1}
          \binom{n-1}{\abs{\pi}_0}
           v(\pi)
    .
  \end{equation}
  Here $\abs{\pi}_0$ is the number of returns to zero of the path
  and $v$ is the valuation \eqref{eq:Flajoletvaluation} from Flajolet's formula.
  Note that the path consisting of $n$ horizontal steps does not contribute.
\end{Theorem}

\begin{proof}
  The proof consists in comparing the expressions for the moments and
  free cumulants in terms of the Boolean cumulants \eqref{eq:booleancumulants},
  which themselves represent the sum over irreducible Motzkin paths
  \eqref{eq:boolean=summotzkinirr}.
  We have then by the multinomial formula
  \begin{equation*}
    \begin{aligned}
      \mu_n
      &= [z^n]M(z) \\
      &= [z^n] \sum_{m=0}^n H(z)^m \\
      &= [z^n] \sum_{m=0}^n \sum_{k_1+k_2+\dots+k_n=m}
                \binom{m}{k_1,\dots,k_n}
                h_1^{k_1}
                \cdots
                h_n^{k_n}
                z^{k_1+2k_2+\dots+nk_n}\\
      &= \sum_{k_1+2k_2+\dots+nk_n=n}
          \binom{k_1+\dots+k_n}%
                {k_1,\dots,k_n}
          h_1^{k_1}
          \cdots
          h_n^{k_n}
    \end{aligned}
    .
  \end{equation*}

The cumulants can be expressed in terms of the moments
with the help of Lagrange's inversion formula \cite{Comtet:1974:advancedcombinatorics}:
  \begin{align*}
    c_n &= -\frac{1}{n-1}
           [z^n] \frac{1}{M(z)^{n-1}}\\
        &= -\frac{1}{n-1}
           [z^n] (1-H(z))^{n-1}\\
        &= -\frac{1}{n-1}
           [z^n] (1 - h_1 z - \dots - h_n z^n)^{n-1}\\
        &= -\frac{1}{n-1}
           [z^n] \sum_{k_0+k_1+\dots+k_n=n-1}
                  \frac{(n-1)!}{k_0!\cdots k_n!}
                  (- h_1)^{k_1}
                  \cdots
                  (- h_n)^{k_n}
                  z^{k_1 + 2 k_2 + \dots + n k_n}\\
        &= -\frac{1}{n-1}
           \sum_{\substack{k_1 + 2 k_2 + \dots + n k_n = n\\ k_1<n}}
            (-1)^{k_1 + \dots + k_n}
            \binom{n-1}{n-1- k_1 - \dots - k_n,k_1,\dots,k_n}
            h_1^{k_1}
            \cdots
            h_n^{k_n}\\
        &= \sum_{\substack{k_1 + 2 k_2 + \dots + n k_n = n\\ k_1<n}}
            \frac{(-1)^{ k_1 + \dots + k_n - 1}}{n-1}
            \binom{n-1}{k_1+\dots+k_n}
            \binom{ k_1 + \dots + k_n}{k_1,\dots,k_n}
            h_1^{k_1}
            \cdots
            h_n^{k_n}\\
        &= \sum_{r=1}^{n-1}
           \sum_{\substack{i_1+\dots+i_r=n\\ i_j\geq 1} }
            \frac{(-1)^{r-1}}{n-1}
            \binom{n-1}{r}
            h_{i_1} h_{i_2} \cdots h_{i_r}
    .
  \end{align*}
  Comparison with the formula for the moments yields the result.
\end{proof}
\begin{Remark}
  Actually by the same argument we can express free cumulants in terms of classical
  cumulants by simply replacing Motzkin paths by \L{}ukasiewicz paths
  with valuation \eqref{eq:Lukasvaluation:classical}.
  However, there are cancellations in the sum.
  There are also lots of cancellations when expressing free cumulants in terms
  of free cumulants themselves using \L{}ukasiewicz paths with valuation 
  \eqref{eq:Lukasvaluation:NC}
  (only one term survives).
  Note however that there are no cancellations in the sum \eqref{eq:freecumulant=motzkin},
  since $\abs{\pi}_0$ is equal to the number of $\alpha_0$'s and $\lambda_1$'s
  and therefore every monomial appears always with the same sign.
\end{Remark}

\section{Hankel determinants}
\label{sec:Hankeldeterminants}

In this section we survey some formulas expressing the Jacobi parameters
$a_n$ and $\lambda_n$ in terms of cumulants.

It is well known (see e.g.\ \cite{Akhiezer:1965:momentproblem})
and can be readily deduced from~\eqref{eq:Dn(x)},
that the Jacobi parameters can be expressed in terms of Hankel determinants
and Hankel minors of the moments as follows:
$$
\lambda_n = \frac{\Delta_{n-2}\Delta_n}{\Delta_{n-1}^2}
$$
and
$$
a_n = \frac{\tilde\Delta_n}{\Delta_n}
      - \frac{\tilde\Delta_{n-1}}{\Delta_{n-1}}
,
$$
where
$$
  \tilde\Delta_n
  = \begin{vmatrix}
      \mu_0    & \mu_1    & \dots & \mu_{n-1} & \mu_{n+1} \\  
      \mu_1    & \mu_2    &       & \mu_{n}   & \mu_{n+2} \\
      \vdots   &          &       & \vdots    & \vdots \\ 
      \mu_n    & \mu_{n+1}& \dots & \mu_{2n-1}& \mu_{2n+1}    
    \end{vmatrix}
.
$$

We are therefore interested to express the minors of the infinite Hankel matrix
$\Hankel_\mu=[\mu_{i+j}]_{i,j\geq0}$ in terms of cumulants.
Namely, given finite sequences of indices
$i_1,i_2,\dots,i_p$ and $j_1,j_2,\dots,j_p$ we will denote
\begin{equation}
  \label{eq:Hankelminor}
  H
  \begin{pmatrix}
    i_1 & i_2 & \dots & i_p\\
    j_1 & j_2 & \dots & j_p
  \end{pmatrix}
  =
  \begin{vmatrix}
    \mu_{i_1+j_1} & \mu_{i_1+j_2} & \dots & \mu_{i_1+j_p} \\
    \mu_{i_2+j_1} & \mu_{i_2+j_2} & \dots & \mu_{i_2+j_p} \\
    \dots\\
    \mu_{i_p+j_1} & \mu_{i_p+j_2} & \dots & \mu_{i_p+j_p}
  \end{vmatrix}
\end{equation}
For example,
$$
\Delta_n
= H\begin{pmatrix}
     0 & 1 & \dots & n\\
     0 & 1 & \dots & n
   \end{pmatrix}
\qquad\qquad
\tilde\Delta_n
= H\begin{pmatrix}
     0 & 1 & \dots & n-1 & n\\
     0 & 1 & \dots & n-1 & n+1
   \end{pmatrix}
.
$$
The main tool is the so-called \emph{Gessel--Viennot theory},
see
\cite{Viennot:1983:polynomes,GesselViennot:1989:determinants,KarlinMcGregor:1950:coincidenceprobabilities,Lindstrom:1973:vectorrepresentations}. 
Let $\Gamma$ be a weighted graph,
that is a graph with vertices $V$ and edges $E$ together with
a valuation $v:E\to R$ on the edges, where $R$ is some (commutative) ring.
The valuation of a path $\omega = (s_0,\dots,s_m)$ is the product
of the valuations of the steps
$v(\omega)=v(s_0,s_1)\,v(s_1,s_2)\,\cdots v(s_{m-1},s_m)$.

Select two sets of distinct vertices $A_i$, $B_i\in V$, $i=1,\dots,n$
and suppose that the set of paths
$$
\Omega_{ij} = \{\omega=(s_0,\dots,s_m) : s_0=A_i, s_m=B_j, v(\omega)\ne 0\}
$$
is finite for every pair of indices $(i, j)$.
We will say that two paths \emph{intersect} if they have a vertex in common.
We will not care about \emph{crossings} of paths, that is,
crossings of edges in a graphical representation.

Define a matrix
$$
a_{ij} = \sum_{\omega\in\Omega_{ij}} v(\omega), \qquad 1\leq i,j \leq n.
$$
Then the determinant of this matrix has the following combinatorial interpretation.
\begin{Proposition}[{\cite{Viennot:1983:polynomes} Prop.~IV.2}]
  \label{thm:gesselviennot}
  \begin{equation}
  \det[a_{ij}]_{n\times n}
  = \sum_{(\sigma; \omega_1,\dots,\omega_n)}
     \sign(\sigma)\, v(\omega_1)\cdots v(\omega_n)
  \end{equation}
  where the sum is over all permutations $\sigma\in \SG_n$, $\omega_i\in\Omega_{i,\sigma_i}$
  and all non-intersecting (but possibly crossing) paths $\omega_i\in\Omega_{i,\sigma(i)}$.
\end{Proposition}

We have seen various examples above where the $n$-th moment is equal
to the sum of the valuations of paths of length $n$.
All our valuations are translation independent and we will interpret
the $(i,j)$ entry of $\Hankel_\mu$ as a sum
$\mu_{i+j} = \sum v(\omega)$ over the set of paths starting
at $(-i,0)$ and ending at $(j,0)$.
The minor \eqref{eq:Hankelminor} is therefore equal to
$$
  H
  \begin{pmatrix}
    i_1 & i_2 & \dots & i_p\\
    j_1 & j_2 & \dots & j_p
  \end{pmatrix}
  = 
  \sum_{\sigma\in\SG_p}
  \sum_{\omega_1,\omega_2,\dots,\omega_p}
   (-1)^{\sign(\sigma)}
   v(\omega_1)\,  v(\omega_2) \cdots v(\omega_p)
$$
where the sum extends over all permutations $\sigma$ of the indices and 
all nonintersecting configurations of paths $\omega_k$ starting at
$(-i_k,0)$ and ending at $(j_{\sigma(k)},0)$.
Because the vertices are ordered,
the sign is $(-1)^{\sign(\pi)}=(-1)^K$, 
where $K$ is the number of crossings of the configuration corresponding
to $\pi$.

For example, for the full Hankel determinant from \eqref{eq:Deltan(x)},
where the moments are interpreted as sum over Motzkin paths as
in Theorem~\ref{thm:Flajolet},
 there is only one non-intersecting configuration whose picture is
(for $n=4$)

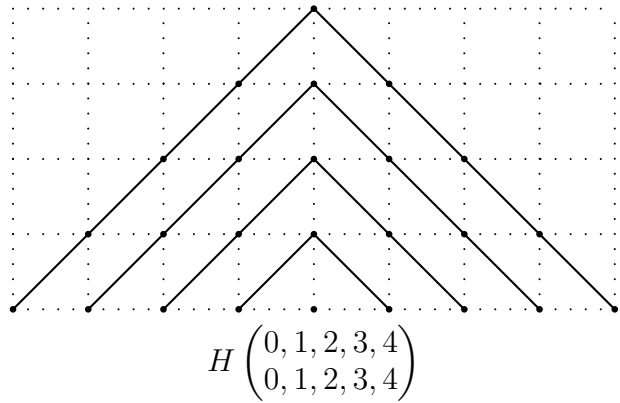
\begin{center}
  \begin{minipage}{10\psunit}
    \begin{pspicture}(-5,0)(5,5)
      \psline(-4, 0)(-3, 1)(-2, 2)(-1, 3)(0, 4)(1, 3)(2, 2)(3, 1)(4, 0)
      \psdots(-4, 0)(-3, 1)(-2, 2)(-1, 3)(0, 4)(1, 3)(2, 2)(3, 1)(4, 0)
      \psline(-3, 0)(-2, 1)(-1, 2)(0, 3)(1, 2)(2, 1)(3, 0)
      \psdots(-3, 0)(-2, 1)(-1, 2)(0, 3)(1, 2)(2, 1)(3, 0)
      \psline(-2, 0)(-1, 1)(0, 2)(1, 1)(2, 0)
      \psdots(-2, 0)(-1, 1)(0, 2)(1, 1)(2, 0)
      \psline(-1, 0)(0, 1)(1, 0)
      \psdots(-1, 0)(0, 1)(1, 0)
      \psdots(0, 0)
      \psgrid[griddots=5,subgriddiv=1,gridlabels=0](-4,0)(4,4)
    \end{pspicture}
    $$
    H\begin{pmatrix}
       0,1,2,3,4\\
       0,1,2,3,4
     \end{pmatrix}
    $$
  \end{minipage}
\end{center}
and therefore
$\Delta_n = \lambda_1^n \lambda_2^{n-1} \cdots \lambda_{n-1}^2 \lambda_n$. 

Note that in this example the point $(0,0)$ is considered as a path of length
zero and is not allowed to lie on any of the other paths.

In contrast, in the case of \L{}ukasiewicz paths with
valuations~\eqref{eq:Lukasvaluation:NC}
and~\eqref{eq:Lukasvaluation:classical}, there may occur crossings of paths
and there are many contributing terms with different signs.
It is also possible that cancellations occur in the sum, e.g.\
$$
v
\left(
  \psset{unit=.5\psunit}
  \begin{minipage}{7\psunit}
  \begin{pspicture}(-3,-1)(4,3)
  \psline(-2, 0)(-1, 1)(0, 2)(1, 1)(2, 0)
  \psdots(-2, 0)(-1, 1)(0, 2)(1, 1)(2, 0)
  \psline(-1, 0)(0, 1)(1, 2)(2, 2)(3, 0)
  \psdots(-1, 0)(0, 1)(1, 2)(2, 2)(3, 0)
  \psline(0, 0)(1, 0)
  \psdots(0, 0)(1, 0)
  \psgrid[griddots=5,subgriddiv=1,gridlabels=0](-2,0)(3,2)
  \end{pspicture}
  \end{minipage}
  \psset{unit=2\psunit}
\right)=-c_1^2 c_2^2 c_3
=
-v
\left(
  \psset{unit=0.5\psunit}
  \begin{minipage}{7\psunit}
  \begin{pspicture}(-3,-1)(4,3)
  \psline(-2, 0)(-1, 1)(0, 2)(1, 2)(2, 0)
  \psdots(-2, 0)(-1, 1)(0, 2)(1, 2)(2, 0)
  \psline(-1, 0)(0, 1)(1, 0)
  \psdots(-1, 0)(0, 1)(1, 0)
  \psline(0, 0)(1, 1)(2, 1)(3, 0)
  \psdots(0, 0)(1, 1)(2, 1)(3, 0)
  \psgrid[griddots=5,subgriddiv=1,gridlabels=0](-2,0)(3,2)
  \end{pspicture}
  \end{minipage}
  \psset{unit=2\psunit}
\right)
$$
are cancelling contributions to
$
H\begin{pmatrix}
   0&1&2\\
   0&1&3
 \end{pmatrix}
$.
Therefore, and because the sums are complicated,
the formulae are of rather limited value.


\begin{thebibliography}{10}

\bibitem{AccardiBozajko:1998:Gaussianization}
Accardi, L. and Bo{\.z}ejko, M., \emph{Interacting {F}ock spaces and
  {G}aussianization of probability measures}, Infin. Dimens. Anal. Quantum
  Probab. Relat. Top. \textbf{1} (1998), 663--670.

\bibitem{Akhiezer:1965:momentproblem}
Akhiezer, N.~I., \emph{The classical moment problem and some related questions
  in analysis}, Hafner Publishing Co., New York, 1965.

\bibitem{AkiyamaYoshida:1999:orthogonal}
Akiyama, M. and Yoshida, H., \emph{The orthogonal polynomials for a linear sum
  of a free family of projections}, Infin. Dimens. Anal. Quantum Probab. Relat.
  Top. \textbf{2} (1999), 627--643.

\bibitem{BozejkoLeinertSpeicher:1996:convolution}
Bo{\.z}ejko, M., Leinert, M. and Speicher, R., \emph{Convolution and limit
  theorems for conditionally free random variables}, Pacific J. Math.
  \textbf{175} (1996), 357--388.

\bibitem{BozejkoSpeicher:1991:psiindependent}
Bo{\.z}ejko, M. and Speicher, R., \emph{$\psi$-independent and symmetrized
  white noises}, Quantum probability \& related topics, World Sci. Publishing,
  River Edge, NJ, 1991, pp.~219--236.

\bibitem{CohenTrenholme:1984:orthogonal}
Cohen, J.~M. and Trenholme, A.~R., \emph{Orthogonal polynomials with a constant
  recursion formula and an application to harmonic analysis}, J. Funct. Anal.
  \textbf{59} (1984), 175--184.

\bibitem{Comtet:1974:advancedcombinatorics}
Comtet, L., \emph{Advanced combinatorics}, D. Reidel Publishing Co., Dordrecht,
  1974.

\bibitem{Cori:1983:wordstrees}
Cori, R., \emph{Words and trees}, Combinatorics on Words (Lothaire, M., ed.),
  Encyclopedia of mathematics, vol.~17, Addison-Wesley Publishing Co., Reading,
  Mass., 1983, pp.~213--227.

\bibitem{Flajolet:1980:continuedfractions}
Flajolet, P., \emph{Combinatorial aspects of continued fractions}, Discrete
  Math. \textbf{32} (1980), 125--161.

\bibitem{FlajoletGuillemin:2000:processes}
Flajolet, P. and Guillemin, F., \emph{The formal theory of birth-and-death
  processes, lattice path combinatorics and continued fractions}, Adv. in Appl.
  Probab. \textbf{32} (2000), 750--778.

\bibitem{Gerl:1984:continuedfraction}
Gerl, P., \emph{Continued fraction methods for random walks on {$\bf N$} and on
  trees}, Probability measures on groups, VII (Oberwolfach, 1983), Lecture
  Notes in Mathematics, vol. 1064, Springer, Berlin, 1984, pp.~131--146.

\bibitem{GesselViennot:1989:determinants}
Gessel, I.~M. and Viennot, X.~G., \emph{Determinants, paths, and plane
  partitions}, Preprint, 1989.

\bibitem{Haagerup:1997:RStransforms}
Haagerup, U., \emph{On {V}oiculescu's ${R}$- and ${S}$-transforms for free
  non-commuting random variables}, Free probability theory (Waterloo, ON,
  1995), Amer. Math. Soc., Providence, RI, 1997, pp.~127--148.

\bibitem{IsmailLetessierMassonValent:1990:processes}
Ismail, M. E.~H., Masson, D.~R., Letessier, J. and Valent, G., \emph{Birth and
  death processes and orthogonal polynomials}, Orthogonal polynomials
  (Columbus, OH, 1989), Kluwer Acad. Publ., Dordrecht, 1990, pp.~229--255.

\bibitem{KarlinMcGregor:1950:coincidenceprobabilities}
Karlin, S. and McGregor, J., \emph{Coincidence probabilities}, Pacific J. Math.
  \textbf{9} (1959), 1141--1164.

\bibitem{Lascoux:2000:motzkin}
Lascoux, A., \emph{Motzkin paths and powers of continued fractions}, S\'em.
  Lothar. Combin. \textbf{44} (2000), Art. B44e, 4 pp.

\bibitem{Lehner:2001:cumulants}
Lehner, F., \emph{Cumulants in noncommutative probability}, in preparation.

\bibitem{Lindstrom:1973:vectorrepresentations}
Lindstr{\"o}m, B., \emph{On the vector representations of induced matroids},
  Bull. London Math. Soc. \textbf{5} (1973), 85--90.

\bibitem{Nica:1995:oneparameter}
Nica, A., \emph{A one-parameter family of transforms, linearizing convolution
  laws for probability distributions}, Comm. Math. Phys. \textbf{168} (1995),
  187--207.

\bibitem{Nica:1996:crossingsembracings}
\bysame, \emph{Crossings and embracings of set-partitions and $q$-analogues of
  the logarithm of the {F}ourier transform}, Discrete Math. \textbf{157}
  (1996), 285--309.

\bibitem{Nica:1996:Rtransforms}
\bysame, \emph{${R}$-transforms of free joint distributions and non-crossing
  partitions}, J. Funct. Anal. \textbf{135} (1996), 271--296.

\bibitem{Promislow:1978:dimension}
Promislow, D., \emph{Dimension of null spaces with applications to group
  rings}, Canad. J. Math. \textbf{30} (1978), 289--300.

\bibitem{Raney:1960:composition}
Raney, G.~N., \emph{Functional composition patterns and power series
  reversion}, Trans. Amer. Math. Soc. \textbf{94} (1960), 441--451.

\bibitem{Speicher:1994:multiplicative}
Speicher, R., \emph{Multiplicative functions on the lattice of noncrossing
  partitions and free convolution}, Math. Ann. \textbf{298} (1994), 611--628.

\bibitem{SpeicherWoroudi:1997:boolean}
Speicher, R. and Woroudi, R., \emph{Boolean convolution}, Free probability
  theory (Waterloo, ON, 1995), Amer. Math. Soc., Providence, RI, 1997,
  pp.~267--279.

\bibitem{Trenholme:1988:Green}
Trenholme, A.~R., \emph{A {G}reen's function for nonhomogeneous random walks on
  free products}, Math. Z. \textbf{199} (1988), 425--441.

\bibitem{Viennot:1983:polynomes}
Viennot, G., \emph{Une th\'eorie combinatoire des polyn\^omes orthogonaux
  g\'en\'eraux}, Notes de conf\'erences, Universit\'e du Qu\'ebec \`a
  Montr\'eal, septembre-octobre 1983.

\bibitem{Viennot:1985:polynomials}
Viennot, G., \emph{A combinatorial theory for general orthogonal polynomials
  with extensions and applications}, Orthogonal polynomials and applications
  (Bar-le-Duc, 1984), Lecture Notes in Mathematics, vol. 1171, Springer,
  Berlin, 1985, pp.~139--157.

\bibitem{VoiculescuDykemaNica:1992:free}
Voiculescu, D.~V., Dykema, K.~J. and Nica, A., \emph{Free random variables},
  CRM Lecture Notes Series, vol.~1, American Mathematical Society, Providence,
  RI, 1992.

\bibitem{Voiculescu:1985:operator}
Voiculescu, D., \emph{Symmetries of some reduced free product ${C}\sp
  \ast$-algebras}, Operator algebras and their connections with topology and
  ergodic theory (Bu\c steni, 1983), Lecture Notes in Mathematics, vol. 1132,
  Springer, Berlin, 1985, pp.~556--588.

\bibitem{Voiculescu:1986:addition}
\bysame, \emph{Addition of certain noncommuting random variables}, J. Funct.
  Anal. \textbf{66} (1986), 323--346.

\bibitem{Voiculescu:2000:lectures}
\bysame, \emph{Lectures on free probability theory}, Lectures on probability
  theory and statistics (Saint-Flour, 1998), Lecture Notes in Math., vol. 1738,
  Springer, Berlin, 2000, pp.~279--349.

\bibitem{vonWaldenfels:1973:approach}
von Waldenfels, W., \emph{An approach to the theory of pressure broadening of
  spectral lines}, Probability and information theory, II, Springer, Berlin,
  1973, pp.~19--69. Lecture Notes in Math., Vol. 296.

\bibitem{vonWaldenfels:1975:intervalpartitions}
\bysame, \emph{Interval partitions and pair interactions}, S\'eminaire de
  Probabilit\'es, IX (Seconde Partie, Univ. Strasbourg, Strasbourg, ann\'ees
  universitaires 1973/1974 et 1974/1975), Springer, Berlin, 1975, pp.~565--588.
  Lecture Notes in Math., Vol. 465.

\bibitem{Zeng:2000:powers}
Zeng, J., \emph{On the powers of {M}otzkin paths}, S\'em. Lothar. Combin.
  \textbf{44} (2000), Art. B44f, 2 pp.

\end{thebibliography}
\bibliographystyle{mamsplain}

\providecommand{\bysame}{\leavevmode\hbox to3em{\hrulefill}\thinspace}

\end{document}